\documentclass{ifacconf}

\usepackage{amsmath, amssymb}
\usepackage{graphicx}      
\usepackage{natbib}        
\usepackage{epstopdf}
\usepackage{overpic}
\usepackage{cancel}
\usepackage{rotating}
\usepackage{url}
\usepackage{caption}
\usepackage{color}
\usepackage{rotating}
\usepackage{multirow}
\usepackage{wrapfig}
\usepackage{mathtools}
\usepackage{arydshln}
\begin{document}
\begin{frontmatter}

\title{Sparse Identification of Nonlinear Dynamics with Control (SINDY\MakeLowercase{c})\thanksref{footnoteinfo}} 

\thanks[footnoteinfo]{ SLB acknowledges support from the U.S. Air Force Center of Excellence on Nature Inspired Flight Technologies and Ideas (FA9550-14-1-0398).  JLP thanks Bill and Melinda Gates for their active support of the Institute of Disease Modeling and their sponsorship through the Global Good Fund.  JNK acknowledges support from the U.S. Air Force Office of Scientific Research (FA9550-09-0174). }

\author[First]{Steven L. Brunton} 
\author[Second]{Joshua L. Proctor}
\author[Third]{J. Nathan Kutz}

\address[First]{Department of Mechanical Engineering, University of Washington, Seattle, WA 98195 USA (e-mail: sbrunton@uw.edu)}
\address[Second]{Institute for Disease Modeling, Bellevue, WA 98004 USA, (e-mail: joproctor@intven.com)}
\address[Third]{Department of Applied Mathematics, University of Washington, Seattle, WA 98195 USA, (e-mail: kutz@uw.edu)}

\begin{abstract}                
Identifying governing equations from data is a critical step in the modeling and control of complex dynamical systems.  
Here, we investigate the data-driven identification of nonlinear dynamical systems with inputs and forcing using regression methods, including sparse regression.  
Specifically, we generalize the sparse identification of nonlinear dynamics (SINDY) algorithm to include external inputs and feedback control.  
This method is demonstrated on examples including the Lotka-Volterra predator--prey model and the Lorenz system with forcing and control.
We also connect the present algorithm with the dynamic mode decomposition (DMD) and Koopman operator theory to provide a broader context.\vspace{-.15in}
\end{abstract}

\end{frontmatter}

\section{Introduction}
The data-driven modeling of complex systems is currently undergoing a revolution.
There is unprecedented availability of high-fidelity measurements from historical records, numerical simulations, and experimental data, and recent developments in machine learning and compressed sensing make it possible to extract more from this data.  
Systems of interest, such as a turbulent fluid, an epidemiological system, a network of neurons, financial markets, or the climate, are high-dimensional, nonlinear, and exhibit multi-scale phenomena in both space and time.  
However, many systems evolve on a low-dimensional attractor that may be characterized by large-scale coherent structures~\citep{HLBR_turb,guckenheimer_holmes}.

System identification comprises a large collection of methods to characterize a dynamical system from data.  
Many techniques in system identification~\citep{ljung:book}, including dynamic mode decomposition (DMD)~\citep{Schmid2008aps,Rowley2009jfm,schmid:2010,Tu2014jcd} and DMD with control (DMDc)~\citep{Proctor2014arxiv}, are designed to handle high-dimensional data with the assumption of linear dynamics; historically, there have been relatively few techniques to identify nonlinear dynamical systems from data.  
However,  DMD has strong connections to nonlinear dynamics through Koopman operator theory~\citep{Koopman1931pnas,Mezic2004physicad,Mezic2005nd}, which spurred significant interest and developments~\citep{Rowley2009jfm,Tu2014jcd,Budivsic2012chaos,Mezic2013arfm}.

A recent breakthrough in nonlinear system identification~\citep{Schmidt2009science} uses genetic programming~\citep{koza1999genetic} to construct families of candidate nonlinear functions for the rate of change of state variables in time.  
A parsimonious model is chosen from this family by finding a Pareto optimal solution that balances model complexity with predictive accuracy.  
In a related modeling framework, we have developed an algorithm for the sparse identification of nonlinear dynamics (SINDY) from data~\citep{Brunton2015arxiv}, relying on the fact that most dynamical systems of interest have relatively few nonlinear terms in the dynamics out of the family of possible terms (i.e., polynomial nonlinearities, etc.). 
This method uses sparsity promoting techniques to find models that automatically balance sparsity in the number of terms with model accuracy.  
An earlier related algorithm~\citep{Wang2011prl} uses compressed sensing~\citep{Donoho2006ieeetit,Candes2006picm,Baraniuk2007ieeespm}, while our algorithm uses sparse regression~\citep{Tibshirani1996lasso} to handle measurement noise and overdetermined cases when we have more time snapshots than state measurements.  

There are many other interesting methods recently developed to incorporate sparsity and nonlinear dynamics~\citep{Schaeffer2013pnas,Ozolicnvs2013pnas,mackey2014compressive}.  
In addition, there are numerous exciting directions in equation-free modeling~\citep{Kevrekidis2003cms}, including the Perron-Frobenius operator~\citep{Froyland2009pd}, cluster reduced-order models based on probabilistic transition between various system behaviors~\citep{Kaiser2014jfm}, and methods for uncertainty quantification and subspace analysis in turbulent flows and the climate~\citep{Majda2007bpnas,Majda2009pnas,Sapsis2013pnas}.

Beyond modeling, a goal for many complex systems is active feedback control, as in many fluid dynamic applications~\citep{Brunton2015amr}. 
Extending the data-driven methods above to disambiguate between the effects of dynamics and actuation is a critical step in developing nonlinear input--output models that are suitable for control design.   
Similar to the extension of dynamic mode decomposition to include the effects of control~\citep{Proctor2014arxiv}, here we extend the SINDY algorithm~\citep{Brunton2015arxiv} to include external inputs and control.  
We also demonstrate the relationship of SINDY, with and without control, to DMD and Koopman methods, concluding that each of these are variations of model identification from data using advanced regression techniques.  

\section{Model identification via regression}
\begin{figure*}
\begin{center}
\begin{overpic}[width=.875\textwidth]{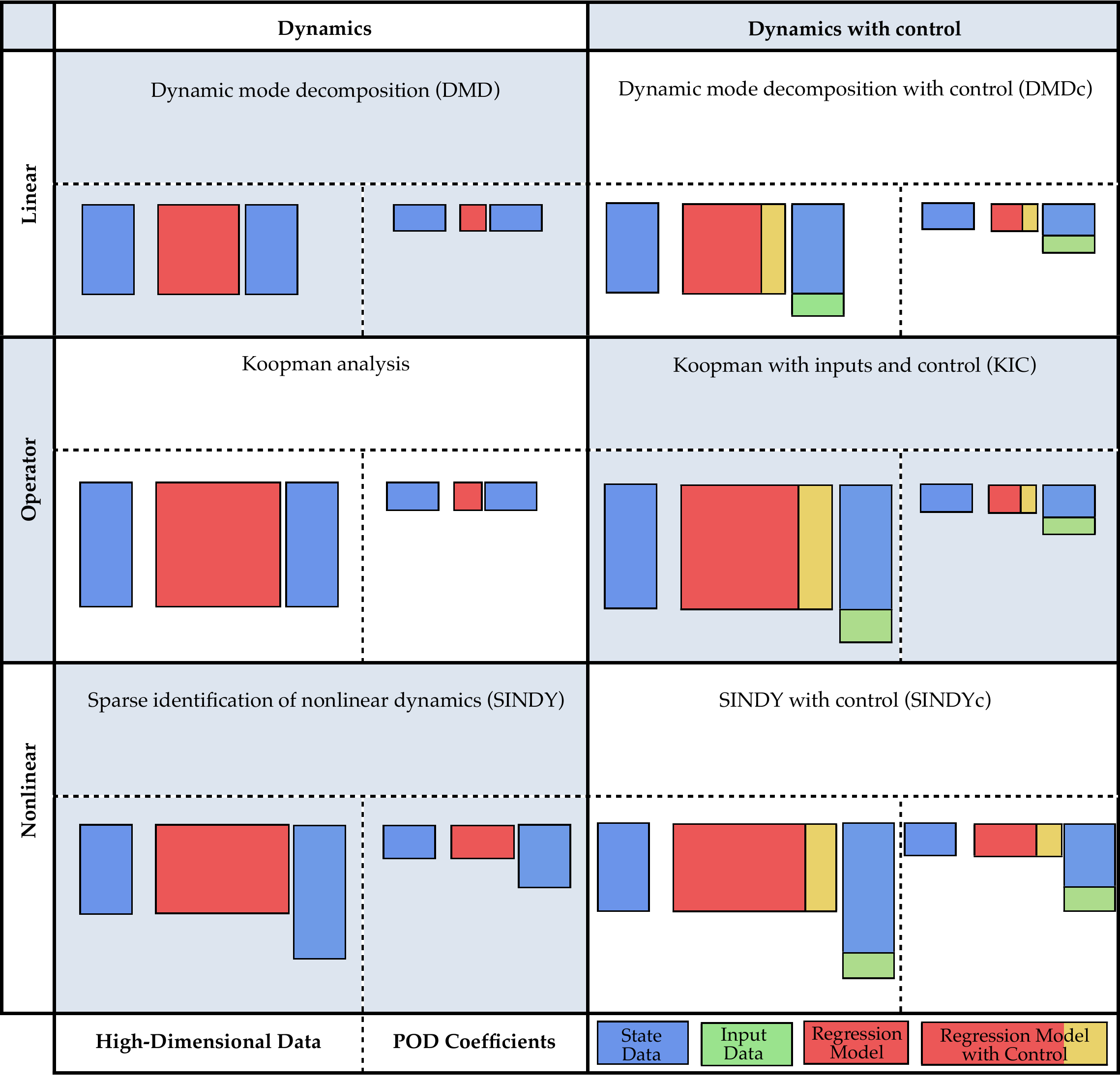}
\put(22,84){$\mathbf{x}_{k+1} = \mathbf{A}\mathbf{x}_k$}
\put(66,84){$\mathbf{x}_{k+1} = \mathbf{A}\mathbf{x}_k+\mathbf{B}\mathbf{u}_k$}
\put(15,59.5){$\mathcal{K}y(\mathbf{x}_k) = y(\mathbf{f}(\mathbf{x}_k)) = y(\mathbf{x}_{k+1})$}
\put(56,59.5){$\mathcal{K}_*y(\mathbf{x}_k,\mathbf{u}_k) = y(\mathbf{f}(\mathbf{x}_k,\mathbf{u}_k),*) = y(\mathbf{x}_{k+1},*)$}
\put(23,29){$\frac{d}{dt}\mathbf{x} = \mathbf{f}(\mathbf{x})$}
\put(69,29){$\frac{d}{dt}\mathbf{x} = \mathbf{f}(\mathbf{x},\mathbf{u})$}
\put(8.5,74){$\mathbf{X}'$}
\put(12.2,74){$=$}
\put(17,74){$\mathbf{A}$}
\put(23.5,74){$\mathbf{X}$}
\put(55.5,74){$\mathbf{X}'$}
\put(59.2,74){$=$}
\put(63.6,74){$\mathbf{A}$}
\put(68.2,74){$\mathbf{B}$}
\put(72.2,74){$\mathbf{X}$}
\put(72.2,68.75){$\mathbf{\Upsilon}$}
\put(39.8,77){\scriptsize$=$}
\put(87.1,77){\scriptsize$=$}
\put(8.5,47.5){$\mathbf{Y}'$}
\put(12.2,47.5){$=$}
\put(18.5,47.5){$\mathbf{K}$}
\put(27,47.5){$\mathbf{Y}$}
\put(55.5,47.5){$\mathbf{Y}'$}
\put(59.,47.5){$=$}
\put(65,47.5){$\mathbf{K}_*$}
\put(76.5,47.5){$\mathbf{Y}$}
\put(76.5,40.){$\mathbf{\Upsilon}$}
\put(39.2,52){\scriptsize$=$}
\put(86.9,51.75){\scriptsize$=$}
\put(8.5,18.5){$\dot{\mathbf{X}}$}
\put(12.2,18.5){$=$}
\put(18.5,18.5){$\mathbf{\Xi}$}
\put(26.2,18.5){\small $\mathbf{\Theta}(\mathbf{X})$}
\put(54.7,18.5){$\dot{\mathbf{X}}$}
\put(58.4,18.5){$=$}
\put(66,18.5){$\mathbf{\Xi}$}
\put(76.7,13.5){\begin{sideways}$\mathbf{\Theta}(\mathbf{X},\mathbf{\Upsilon})$\end{sideways}}
\put(38.9,21){\scriptsize$=$}
\put(85.5,21.25){\scriptsize$=$}
\end{overpic}
\vspace{-.1in}
\caption{\small Overview of various methods that use regression to identify dynamics from data.}\label{fig:overview}
\end{center}
\end{figure*}

Here we review various techniques in system identification, including dynamic mode decomposition (DMD), Koopman analysis, and the sparse identification of nonlinear dynamics (SINDY).  Each of these methods is cast as a regression problem of data onto models, and the schematic overview of these methods is shown in Fig.~\ref{fig:overview}.

\subsection{Dynamic mode decomposition}
The dynamic mode decomposition (DMD) originated in the fluids community to extract spatial-temporal coherent structures from fluid data sets~\citep{Schmid2008aps,Rowley2009jfm,schmid:2010,Tu2014jcd}.    
DMD modes are spatially coherent and oscillate at a fixed frequency and/or growth or decay rate.  
Since fluids data is typically high-dimensional, DMD is built on the proper orthogonal decomposition (POD)~\citep{HLBR_turb}, effectively recombining POD modes in a linear combination to enforce the temporal coherence.  

First, we collect multiple snapshots of high-dimensional fluid data in time $\mathbf{x}_k = \mathbf{x}(k\Delta t)\in\mathbb{R}^n$, where $n$ represents the number of spatial measurements, which may easily represent millions or billions of degrees of freedom.  In DMD, we seek a linear operator $\mathbf{A}$ that approximately relates these snapshots, at least for short periods of time:
\begin{eqnarray}
\mathbf{x}_{k+1} &\approx& \mathbf{A}\mathbf{x}_k.\label{Eq:LinearDynamics}
\end{eqnarray}

If we collect $m+1$ snapshots and arrange in two matrices:
\begin{eqnarray}
\mathbf{X}= \begin{bmatrix} \vline & \vline & & \vline \\
\mathbf{x}_1 & \mathbf{x}_2 & \cdots & \mathbf{x}_m\\
\vline & \vline & & \vline
\end{bmatrix} ,~~ \mathbf{X}'= \begin{bmatrix} \vline & \vline & & \vline \\
\mathbf{x}_2 & \mathbf{x}_3 & \cdots & \mathbf{x}_{m+1}\\
\vline & \vline & & \vline
\end{bmatrix},\label{Eq:DMDStates}
\end{eqnarray}
it is possible to related these matrices by:
\begin{eqnarray}
\mathbf{X}' \approx \mathbf{A}\mathbf{X}.\label{Eq:DMD1}
\end{eqnarray}

In principle, for low-dimensional data it is possible to solve directly for the best-fit linear operator $\bar{\mathbf{A}}$ that minimizes $\|\mathbf{X}'-\bar{\mathbf{A}}\mathbf{X}\|_{F}$ using a least-squares regression, where $\|\cdot\|_F$ is the Frobenius norm.  Numerically, the singular value decomposition (SVD) is used to apply the pseudo-inverse of $\mathbf{X}$ to both sides of Eq.~\eqref{Eq:DMD1}.  
However, when the state dimension $n$ is large, then $\mathbf{A}$ is high-dimensional with $n^2$ elements, and might not be representable computationally.  
Instead, we apply the proper orthogonal decomposition to the data $\mathbf{X}$ and compute a reduced operator $\tilde{\mathbf{A}}$ that acts on POD coefficients.  
It is possible to reconstruct the leading eigenvalues and eigenvectors of the high-dimensional $\mathbf{A}$ matrix from the eigendecomposition of $\tilde{\mathbf{A}}$.  

\begin{enumerate}
\item Compute the economy-sized SVD of $\mathbf{X}$:
\begin{eqnarray}
\mathbf{X} = \mathbf{U}\mathbf{\Sigma}\mathbf{V}^*,
\end{eqnarray}
where $\mathbf{U}\in\mathbb{R}^{n\times m}$, $\mathbf{\Sigma}\in\mathbb{R}^{m\times m}$, and $\mathbf{V}\in\mathbb{R}^{m\times m}$.
\item Compute the projection of the least-square solution $\bar{\mathbf{A}} = \mathbf{X}'\mathbf{X}^{\dagger}$ onto POD modes, given by the columns of $\mathbf{U}$, where $\mathbf{X}^{\dagger}=\mathbf{V}\mathbf{\Sigma}^{-1}\mathbf{U}^*$ is the psuedo-inverse:
\begin{eqnarray}
\tilde{\mathbf{A}} = \mathbf{U}^*\bar{\mathbf{A}}\mathbf{U} = \mathbf{U}^*\mathbf{X}'\mathbf{V}\mathbf{\Sigma}^{-1}.
\end{eqnarray}
Note that $\tilde{\mathbf{A}}$ is an $m\times m$ matrix, where $m$ is the number of time snapshots; this matrix advances POD coefficients forward in time.
\item Compute the eigendecomposition of $\tilde{\mathbf{A}}$:
\begin{eqnarray}
\tilde{\mathbf{A}}\mathbf{W} = \mathbf{W}\mathbf{\Lambda}.
\end{eqnarray}
\item The eigenvalues in $\mathbf{\Lambda}$ are also eigenvalues of the full $\mathbf{A}$ matrix, and these are called \emph{DMD eigenvalues}.  The corresponding eigenvectors of $\mathbf{A}$, called \emph{DMD modes}, are constructed as~\citep{Tu2014jcd}:
\begin{eqnarray}
\mathbf{\Phi} = \mathbf{X}'\mathbf{V}\mathbf{\Sigma}^{-1}\mathbf{W}.
\end{eqnarray}
\end{enumerate}
It is also possible to truncate the SVD at order $r$, retaining only the first $r$ POD modes, and resulting in an $r\times r$ matrix $\tilde{\mathbf{A}}$.  
The DMD modes $\boldsymbol{\phi}$ are spatially coherent and oscillate and/or grow or decay at the fixed frequency $\lambda$.  

The dynamic mode decomposition has been applied to a wide range of problems including fluid mechanics~\citep{Rowley2009jfm,Tu2014jcd,Mezic2013arfm}, epidemiology~\citep{Proctor2015ih}, neuroscience~\citep{brunton2014extracting}, robotics~\citep{Berger2014ieee}, and video processing~\citep{grosek2014,Erichson2015arxiv}.  
However, many of these applications have the ultimate goal of closed-loop feedback control.

\subsection{Dynamic mode decomposition with control}

To disambiguate the effect of internal dynamics from actuation or external inputs, the dynamic mode decomposition with control (DMDc) was developed~\citep{Proctor2014arxiv}.
In DMDc, the linear state dynamics in Eq.~\eqref{Eq:LinearDynamics} are augmented to include the effect of actuation inputs $\mathbf{u}$:
\begin{eqnarray}
\mathbf{x}_{k+1} &\approx& \mathbf{A}\mathbf{x}_k+\mathbf{B}\mathbf{u}_k.\label{Eq:LinearDynamicsWithControl}
\end{eqnarray}
We still collect the state snapshots from Eq.~\eqref{Eq:DMDStates}, but now we collect an additional matrix for the control history:
\begin{eqnarray}
\mathbf{\Upsilon} = \begin{bmatrix} \vline & \vline & & \vline \\
\mathbf{u}_1 & \mathbf{u}_2 & \cdots & \mathbf{u}_m\\
\vline & \vline && \vline \end{bmatrix}.
\end{eqnarray}

In DMDc, $\mathbf{A}$ and $\mathbf{B}$ are approximated from data via:
\begin{eqnarray}
\mathbf{X}' \approx \begin{bmatrix} \mathbf{A} & \mathbf{B}\end{bmatrix} \begin{bmatrix} \mathbf{X}\\ \mathbf{\Upsilon}\end{bmatrix}.
\end{eqnarray}

\subsection{Koopman analysis}
DMD is connected to nonlinear systems via the Koopman operator~\citep{Mezic2004physicad,Mezic2005nd,Rowley2009jfm,Tu2014jcd}.  
The Koopman operator~\citep{Koopman1931pnas} is an infinite-dimensional \emph{linear} operator that describes how a measurement function $y(\mathbf{x})$ evolves through \emph{nonlinear} dynamics:
\begin{eqnarray}
\mathbf{x}_{k+1} &=& {\bf f}(\mathbf{x}_k).
 \end{eqnarray}
 The Koopman operator $\mathcal{K}$ acts on the Hilbert space of scalar measurement functions $y(\mathbf{x})$ as:
 \begin{eqnarray}
\mathcal{K} y(\mathbf{x}_k) =y(\mathbf{f}(\mathbf{x}_k)) = y(\mathbf{x}_{k+1}). 
\end{eqnarray}
That is, the Koopman operator acts on $y$ by the composition of $y$ with the dynamic update $\mathbf{f}$.  

The DMD algorithm approximates the spectrum of the Koopman operator using linear observable functions (i.e., the observable functions are linear functions of the state, as in $y(\mathbf{x}_k) = \mathbf{x}_k$).  However, it was recently shown that linear measurements are not sufficiently rich to analyze nonlinear systems~\citep{Williams2014arxivA}, resulting in the extended DMD (eDMD), which performs a similar DMD regression, but on an augmented data matrix including nonlinear state measurements.  Since this algorithm is expensive numerically, a kernel trick was implemented to make the eDMD method as computationally efficient as standard DMD~\citep{Williams2015jnls}.

\subsection{Koopman with inputs and control}
Similar to how DMD was extended to include inputs and control, Koopman analysis has recently been extended to include inputs and control~\citep{Proctor2016arxiv}.  In this Koopman with inputs and control (KIC) framework, scalar measurements of the state and control $y(\mathbf{x},\mathbf{u})$ are advanced through nonlinear dynamics with control:
\begin{eqnarray}
\mathbf{x}_{k+1} = \mathbf{f}(\mathbf{x}_k,\mathbf{u}_k).
\end{eqnarray}
The Koopman with control operator $\mathcal{K}_*$ is given by:
\begin{eqnarray}
\mathcal{K}_* y(\mathbf{x}_k,\mathbf{u}_k) =y(\mathbf{f}(\mathbf{x}_k,\mathbf{u_k}),*) = y(\mathbf{x}_{k+1},*). 
\end{eqnarray}
It is important to note that there is a parameterized family of Koopman with control operators $\mathcal{K}_*$, as there is a choice of which future control input $*$ to use.  It has been shown that Koopman with inputs and control reduces to DMDc for linear dynamical systems, much as Koopman analysis is numerically computed using DMD for linear systems.

\subsection{Sparse identification of nonlinear dynamics (SINDY)}
The SINDY algorithm identifies fully nonlinear dynamical systems from measurement data.  
This relies on the fact that many dynamical systems have relatively few terms in the right hand side of the governing equations:
\begin{eqnarray}
\frac{d}{dt}\mathbf{x} = \mathbf{f}(\mathbf{x}).\label{Eq:NonlinearDynamics}
\end{eqnarray}
Given a library of candidate nonlinear functions, 
\begin{eqnarray}
\mathbf{\Theta}^T(\mathbf{X}) = \begin{bmatrix}
\rule[2.2pt]{4em}{0.4pt} & \mathbf{1} & \rule[2.2pt]{4em}{0.4pt} \\ 
\rule[2.2pt]{4em}{0.4pt} & \mathbf{X} & \rule[2.2pt]{4em}{0.4pt}\\
\rule[2.2pt]{4em}{0.4pt} & \mathbf{X}^2 & \rule[2.2pt]{4em}{0.4pt}\\
 & \vdots &\\
 \rule[2.2pt]{4em}{0.4pt} & \sin(\mathbf{X}) & \rule[2.2pt]{4em}{0.4pt}\\
 \rule[2.2pt]{4em}{0.4pt} & \sin(2\mathbf{X}) & \rule[2.2pt]{4em}{0.4pt}\\
 & \vdots &
\end{bmatrix},
\end{eqnarray}
where $\mathbf{X}$ is the same data matrix as in Eq.~\eqref{Eq:DMDStates}, we may write our dynamical system as:
\begin{eqnarray}
\dot{\mathbf{X}} = \mathbf{\Xi}\mathbf{\Theta}^T(\mathbf{X}).\label{Eq:SINDY}
\end{eqnarray}
The coefficients $\mathbf{\Xi}$ in this library are \emph{sparse} for most dynamical systems.  
Therefore, we employ sparse regression to identify a sparse $\mathbf{\Xi}$ corresponding to the fewest nonlinearities in our library that give good model performance.  
Choosing a library of candidate dynamics is a crucial choice int he SINDY algorithm.  The algorithm may be extended to include support for more general nonlinearities.  It may also be possible to test different libraries (polynomials, trigonometric functions, etc.) and also incorporate partial knowledge of the physics (fluids vs. quantum mechanics, etc.).

Notice that if $\mathbf{\Theta}^T(\mathbf{X})=\mathbf{X}$, then Eq.~\eqref{Eq:SINDY} is equivalent to DMD with $\mathbf{\Xi}=\mathbf{A}$.
Each row of Eq.~\eqref{Eq:SINDY} represents a row in Eq.~\eqref{Eq:NonlinearDynamics}, and the sparse vector of coefficients $\boldsymbol{\xi}_k$ corresponding to the $k$-th row of $\mathbf{\Xi}$ is found using a sparse regression algorithm, such as LASSO~\citep{Tibshirani1996lasso}:
\begin{eqnarray}
\boldsymbol{\xi}_k = \text{argmin}_{\boldsymbol{\xi}_k}\|\dot{\mathbf{X}}_k-\boldsymbol{\xi}_k\mathbf{\Theta}^T(\mathbf{X})\|_2 + \alpha \|\boldsymbol{\xi}_k\|_1,
\end{eqnarray}
where $\dot{\mathbf{X}}_k$ represents the $k$-th row of $\dot{\mathbf{X}}$.  The $\|\cdot\|_1$ term promotes sparsity in the coefficient vector $\boldsymbol{\xi}_k$.  The parameter $\alpha$ is selected to identify the Pareto optimal model that best balances low model complexity with accuracy.  A coarse sweep of $\alpha$ is performed to identify the rough order of magnitude where terms are eliminated and where error begins to increase.  Then this parameter sweep may be refined.  

To approximate derivatives from noisy state measurements, the SINDY algorithm uses the total variation regularized derivative~\citep{Rudin1992physd,Chartrand2011isrnam}.

\section{Sparse identification of nonlinear dynamics with control (SINDY\MakeLowercase{c})}
Here, we generalize the SINDY method to include inputs and control.  In particular, we now consider the nonlinear dynamical system with inputs $\mathbf{u}$:
\begin{eqnarray}
\frac{d}{dt}\mathbf{x} &=& {\bf f}(\mathbf{x},\mathbf{u}).\label{Eq:NonlinearDynamicsWithControl}
 \end{eqnarray}
 
 The SINDY algorithm is readily generalized to include actuation, as this merely requires building a larger library $ \mathbf{\Theta}(\mathbf{x},\mathbf{u})$ of candidate functions that include $\mathbf{u}$; these functions can include nonlinear cross terms in $\mathbf{x}$ and $\mathbf{u}$.  
 This extension requires measurements of the state $\mathbf{x}$ as well as the input signal $\mathbf{u}$.  
 This generalization is shown in Fig.~\ref{fig:overview} in terms of the overarching regression framework.  
 
If the signal $\mathbf{u}$ corresponds to an external forcing, then we solve for the sparse coefficients $\mathbf{\Xi}$ in the following:
\begin{eqnarray}
\dot{\mathbf{X}} = \mathbf{\Xi}\mathbf{\Theta}^T(\mathbf{X},\boldsymbol{\Upsilon}).\label{Eq:SINDYc}
\end{eqnarray}
However, if the signal $\mathbf{u}$ corresponds to a feedback control signal, so that $\mathbf{u} = \mathbf{k}(\mathbf{x})$, then it is impossible to disambiguate the effect of the feedback control $\mathbf{u}$ with internal feedback terms $\mathbf{k}(\mathbf{x})$ within the dynamical system; namely, the SINDY regression becomes ill-conditioned.  In this case, we may identify the actuation $\mathbf{u}$ as a function of the state:
\begin{eqnarray}
\boldsymbol{\Upsilon} = \mathbf{\Xi}_u\mathbf{\Theta}^T(\mathbf{X}).
\end{eqnarray}  
To identify the coefficients $\mathbf{\Xi}$ in Eq.~\eqref{Eq:SINDYc}, we perturb the signal $\mathbf{u}$ to allow it to be distinguished from $\mathbf{k}(\mathbf{x})$ terms.  
This may be done by injecting a sufficiently large white noise signal, or occasionally kicking the system with a large impulse or step in $\mathbf{u}$.  
An interesting future direction would be to design input signals that \emph{aid} in the identification of the dynamical system in Eq.~\eqref{Eq:NonlinearDynamicsWithControl} by perturbing the system in directions that yield high-value information.

\section{Example systems}
Here, we demonstrate the SINDY with control algorithm on a simple example predator-prey model with forcing and on the Lorenz equations with external forcing and control.  

\subsection{Predator-prey model}
A predator-prey model with forcing is given by:  
\begin{subequations}
\begin{eqnarray}
\dot x_1 & = & a x_1 - b x_1 x_2 + u^2,\\
\dot x_2 & = & -c x_2 + d x_1x_2 
\end{eqnarray}\label{Eq:LV}
\end{subequations}
The variable $x_1$ represents the size of the prey population and $x_2$ represents the size of the predator population; the prey species is actuated with $u^2$.  The parameters $a,b,c,$ and $d$ represents the various growth/death rates, the effect of predation on the prey population, and the growth of predators based on the size of the prey population.

\begin{figure}
\begin{center}
\vspace{-.1in}
\includegraphics[width=.415\textwidth]{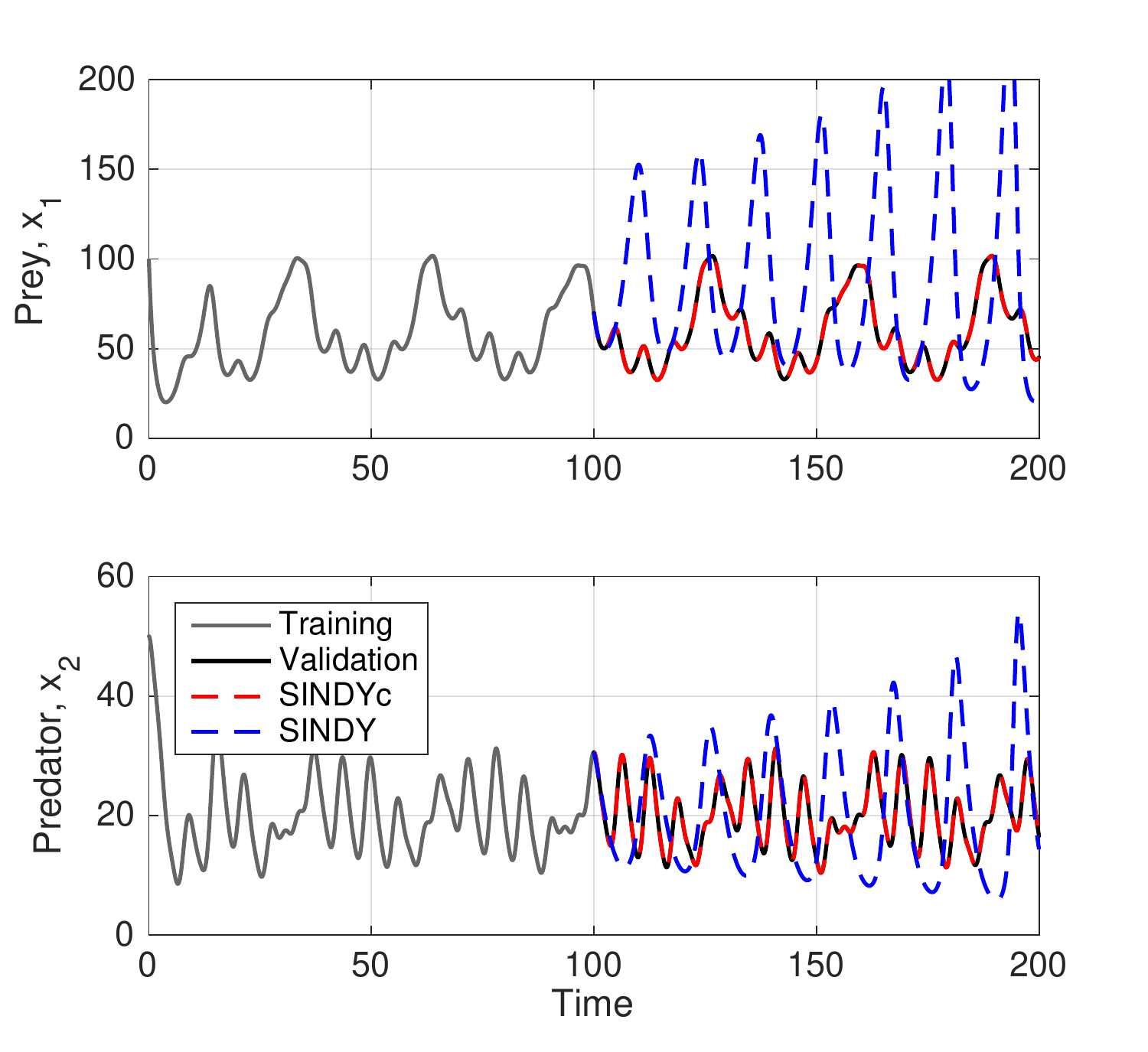}
\end{center}
\vspace{-.175in}
\caption{\small SINDY and SINDYc model predictions for the force Lotka-Volterra system in Eq.~\eqref{Eq:LV}.  The training data consists of the Lotka-Volterra system with periodic forcing.}\label{fig:LotkaRecon}
\end{figure}

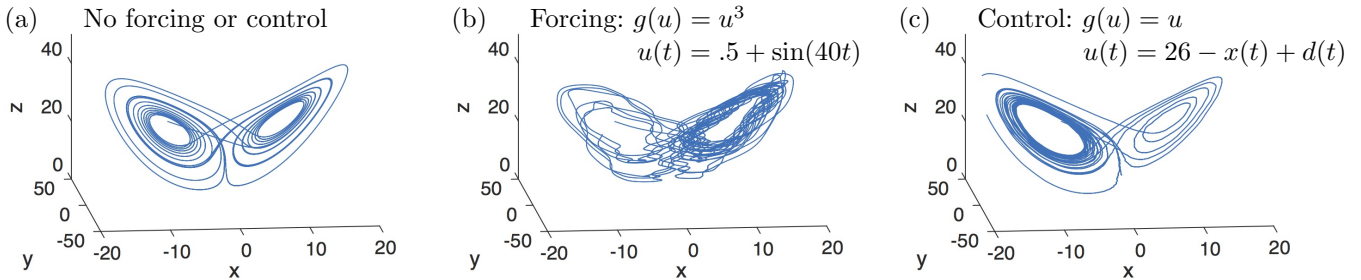
\begin{figure*}
\begin{center}
\begin{tabular}{ccc}
\begin{overpic}[width=.3\textwidth]{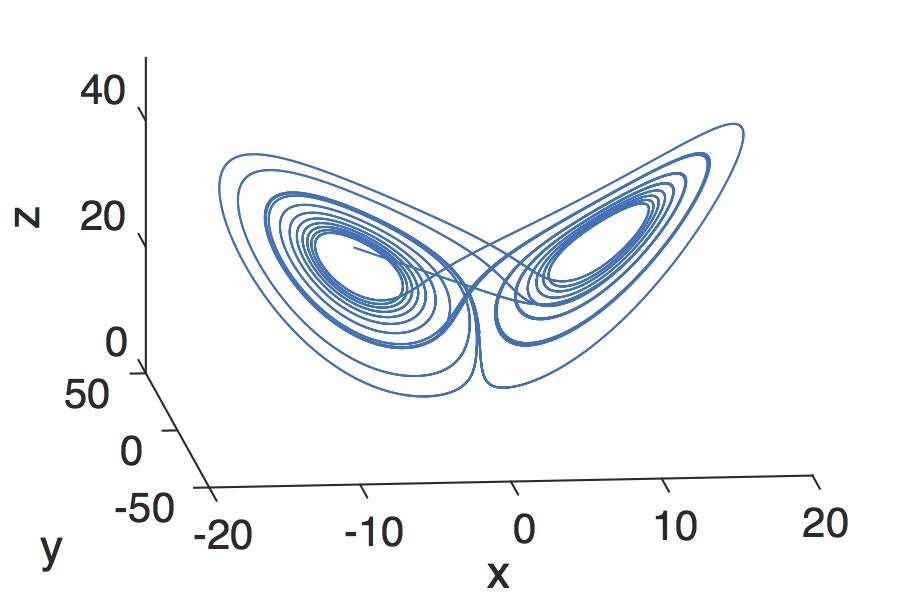}
\put(0,62){(a)}
\put(19,62){No forcing or control}
\end{overpic}
&
\begin{overpic}[width=.3\textwidth]{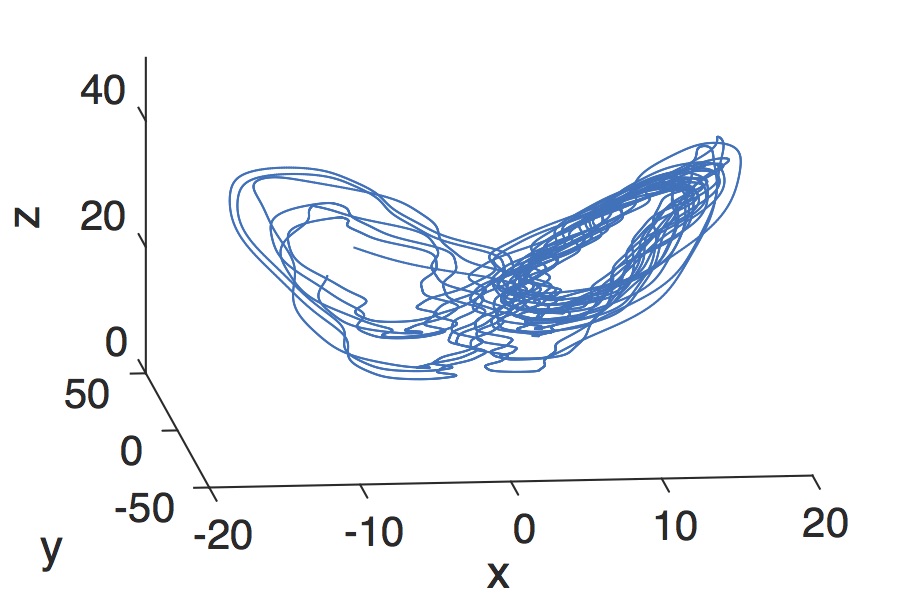}
\put(0,62){(b)}
\put(19,62){Forcing: $g(u)=u^3$}
\put(44.7,54){$u(t) = .5+\sin(40t)$}
\end{overpic}
&
\begin{overpic}[width=.3\textwidth]{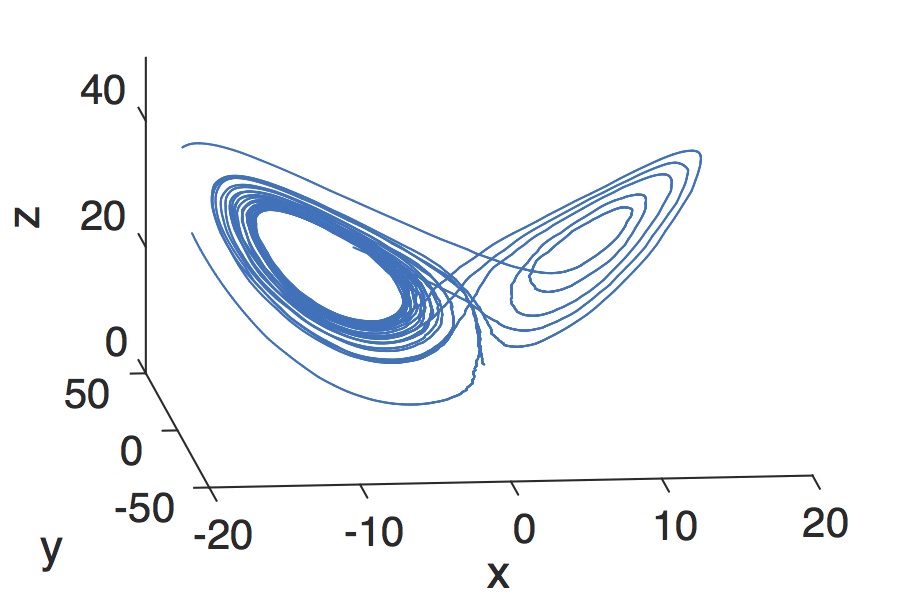}
\put(0,62){(c)}
\put(19,62){Control: $g(u)=u$}
\put(44.5,54){$u(t) = 26-x(t)+d(t)$}
\end{overpic}
\end{tabular}
\end{center}
\vspace{-.1in}
\caption{\small Lorenz system without forcing or control (a), with forcing (b), and with feedback control (c).  In each case, the system is integrated for $50$ time units with a $\Delta t=.001$ with the parameters $\sigma=10$, $\beta=8/3$, and $\rho=28$.}\label{fig:Lorenz}
\end{figure*}

In this example, we force the system sinusoidally with $u(t) = 2\sin(t)+2\sin(t/10)$, and the population response is shown in Fig.~\ref{fig:LotkaRecon} (grey and black).  
The first $100$ time units are used to train the SINDY and SINDYc algorithms, after which they are validated on the next $100$ time units of forced data.  The naive application of SINDY without knowledge of the input results in an unstable model (blue), while the SINDYc algorithm correctly identifies the model structure and parameters in Eq.~\eqref{Eq:LV} to within machine precision in the absence of measurement noise; the SINDYc reconstruction is shown in red.  

\subsection{Lorenz equations}
We also test the SINDYc method on the Lorenz equations:
\begin{subequations}
\begin{eqnarray}
\dot x & = & \sigma(y-x) + g(u) \\
\dot y & = & x(\rho-z) -y \\
\dot z & = & xy-\beta z.
\end{eqnarray}\label{Eq:Lorenz}
\end{subequations}
These equations are examined with various forcing and control models, as shown in Fig.~\ref{fig:Lorenz}.  
In the case of an external forcing, as in Fig.~\ref{fig:Lorenz} (b), the SINDYc algorithm correctly identifies the model and nonlinear input terms.  

In the case that the Lorenz system is being actively controlled by state feedback, as in Fig.~\ref{fig:Lorenz} (c), we must add a perturbation signal $d(t)$ to the input to disambiguate the effect of state feedback via $u$ from internal dynamics.  For this problem, we use an additive white noise process.  
In this example, we train the models using $20$ time units of controlled data, and validate them on another $20$ time units where we switch the forcing to a periodic signal $u(t) = 50\sin(10t)$.  The SINDY algorithm does not capture the effect of actuation, while SINDYc correctly identifies the model and predicts the behavior in response to a new forcing that was not used in the training data.  

\begin{figure}
\begin{center}
\vspace{-.15in}
\includegraphics[width=.425\textwidth]{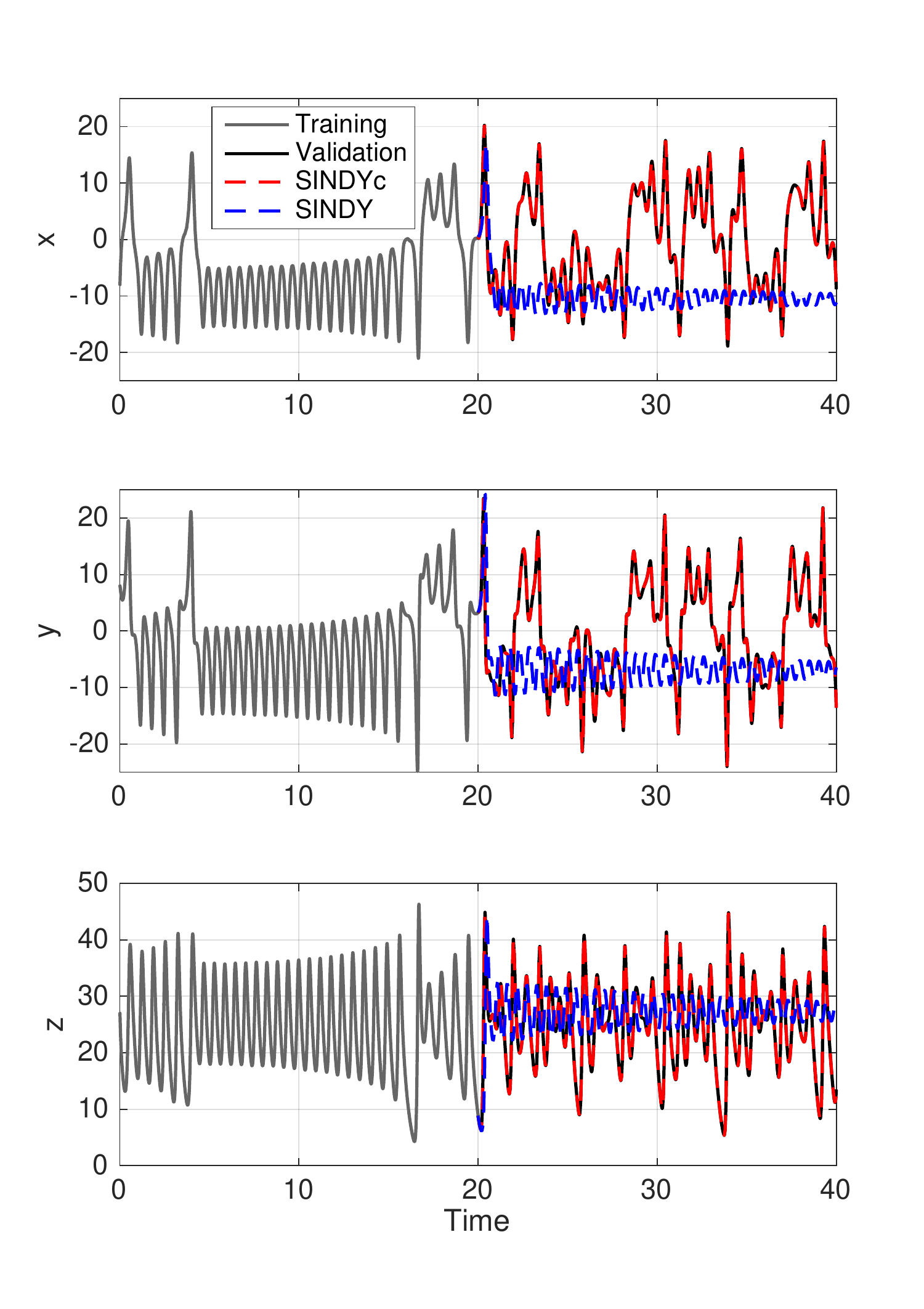}
\end{center}
\vspace{-.3in}
\caption{\small SINDY and SINDYc  predictions for the controlled Lorenz system in Eq.~\eqref{Eq:Lorenz}.  Training data consists of the Lorenz system with state feedback as in Fig.~\ref{fig:Lorenz} (c).  After the training period, the input $u$ switches to a periodic signal $u(t) = 50\sin(10t)$.}\label{fig:LorenzRecon}
\end{figure}

\section{Discussion}
In this work, we have generalized the sparse identification of nonlinear dynamics (SINDY) algorithm to include inputs and control.  
This involved generalizing the library of candidate nonlinear terms to include functions not only of the state $\mathbf{x}$, but also of the input $\mathbf{u}$, including cross terms between state and input.  
This new algorithm is cast in an overarching regression framework in Fig.~\ref{fig:overview}, relating it to other algorithms that determine models from data, including dynamics mode decomposition (DMD), DMD with control, extended DMD, and Koopman analysis.  

The new method has been tested on a predator prey model and the Lorenz system with various forcing and control models.  
The proposed algorithm should scale to the same class of problems where SINDY is useful, since they are built on the same computational architecture.

There are a number of interesting directions to extend this work. 
First, it is important to determine optimal strategies to disambiguate the effect of a state-feedback control signal from internal state dynamics; this may be achieved by additive white noise on the input signal or occasional kicks to the system, but understanding the tradeoffs and benefits of these strategies will be useful.  
More importantly, it is likely possible to design input sequences that optimally probe complex systems to extract high-value information that will be useful to characterize the system.  
For example, perhaps perturbing some systems off-attractor will provide valuable information about nonlinear terms in the dynamics if the on-attractor data may strongly resemble a linear system.  
If the state and control variables have different levels of sparsity, it may be possible to use a weighted convex optimization to penalize the state and control sparsity separately.  
It may also important to improve the model identification if the control law $\mathbf{u} = \mathbf{k}(\mathbf{x})$ is known.  
These are promising areas of current and future research.

\small
\bibliography{nolcos}         
\end{document}